\documentclass[final,5p,times,twocolumn]{elsarticle}

\usepackage{amsfonts, amssymb, amsbsy, amsmath, mathrsfs, amsthm}
\usepackage{tikz}
\usepackage{hyperref}
\usepackage{svg}
\usetikzlibrary{arrows.meta,arrows}
\usepackage{setspace}
\usepackage{mathabx}
\usepackage{algorithm}
\usepackage{algorithmicx}
\usetikzlibrary{shapes}

\algdef{SE}[SUBALG]{Indent}{EndIndent}{}{\algorithmicend\ }%
\algtext*{Indent}
\algtext*{EndIndent}

\theoremstyle{plain}
\newtheorem{theorem}{Theorem}[section]

\newtheorem{proposition}[theorem]{Proposition}

\newtheorem{definition}[theorem]{Definition}

\makeatletter
\def\ps@pprintTitle{%
 \let\@oddhead\@empty
 \let\@evenhead\@empty
 \def\@oddfoot{}%
 \let\@evenfoot\@oddfoot}
\makeatother

\journal{Operations Research Letters}

\begin{document}

\begin{frontmatter}


\title{Risk-averse decision strategies for influence diagrams using rooted junction trees}


\author[label1]{Olli Herrala}
\author[label1]{Topias Terho}
\author[label1]{Fabricio Oliveira\corref{cor1}}
\cortext[cor1]{Corresponding author: fabricio.oliveira@aalto.fi}
\affiliation[label1]{organization={Department of Mathematics and Systems Analysis},
            addressline={Aalto University, School of Science}, 
            postcode={FI-00076 Aalto}, 
            country={Finland}}

\begin{abstract}
This paper presents how a mixed-integer programming (MIP) formulation for influence diagrams, based on a gradual rooted junction tree representation of the diagram, can be generalized to incorporate risk considerations such as conditional value-at-risk and chance constraints. We present two algorithms on how targeted modifications can be made to the underlying influence diagram or to the gradual rooted junction tree representation to enable our reformulations. We present computational results comparing our reformulation with another MIP formulation for influence diagrams.
\end{abstract}




\end{frontmatter}

\section{Introduction}
\label{sec:intro}


An influence diagram (ID) \citep{howard2005influence} is an intuitive structural representation of a decision problem with uncertainties and interdependencies between random events, decisions and consequences. Traditional solution methods for IDs \citep{shachter2010solving} often require strong assumptions such as the no-forgetting assumption. \citet{lauritzen2001representing} present the notion of a limited memory influence diagram (LIMID) that, albeit more general in terms of representation capabilities, does not satisfy the no-forgetting assumption and, therefore, is not amenable to these traditional methods. 

The algorithms presented in the literature for solving decision problems represented as IDs are mostly suited only to problems where an expected utility function is maximized and risk is not explicitly constrained. Thus, often risk considerations are encoded in the utility function itself, by making it concave using, e.g., utility extraction techniques \citep{davidson1957decision,geissel2018optimal,  nielsen2004learning}.
Utility functions often represent monetary values, such as costs or revenues. In that case, maximizing expected utility assumes that the decision-maker has a risk-neutral stance. However, decision-makers may still have different risk tolerance profiles, which must be represented in the decision process.

There are numerous ways to incorporate risk aversion into decision models without requiring utility extraction techniques. A typical method is to minimize a risk measure instead of expected utility \citep{HomemDeMelloPagnocelli2016}. A commonly used measure is the \textit{conditional Value-at-Risk} (CVaR), which measures the expected loss value in the $\alpha$-tail, $\alpha$ being a confidence level parameter \citep{RockafellerUryasev2000}. Another typical way of incorporating risk aversion is to use constraints such as those related to chance events or budget violations \citep{CharnesCooper1959}. Both mentioned methods have been used widely in various applications (See, e.g., \citep{FilippiEtAl2020,KhassibaEtAl2020,XuEtAl2017}). Directly optimizing a risk measure within an ID is challenging as, unlike expected utility, it prevents the use of methods that construct the optimal strategy by computing locally optimal strategies at individual decision nodes.

Recently, two different mixed-integer programming (MIP) reformulations for IDs have emerged, likely stemming from the considerable computational improvements in MIP solution methods. The reformulation considered in this paper is originally presented by \citet{parmentier2020integer}, where the authors first show how to convert a LIMID representing an expected utility maximization problem into a gradual rooted junction tree. This junction tree consists of clusters of nodes from the LIMID and is reformulated as a MIP problem using marginal probability distributions of nodes within each cluster. However, \citet{parmentier2020integer} only consider expected utility maximization and do not show how risk can be accounted for in their formulation.

In contrast, \citet{salo2022decision} present decision programming, which reformulates a LIMID as a mixed-integer linear programming (MILP) formulation without the intermediate clustering step of forming a junction tree. The decision programming formulation used in this paper is the one presented in \citep{hankimaa2023decisionprogrammingjl}, which improves that originally proposed in \citep{salo2022decision} by means of valid inequalities and reformulations.

In the context of MIP formulations for influence diagrams, the main advantage of decision programming is that its formulation can be adapted to minimize risk measures, including CVaR. Similarly, chance-, logical, and budget constraints can easily be incorporated into their MILP formulation, as discussed by \citet{hankimaa2023decisionprogrammingjl}. However, flexibility comes at a cost of computational efficiency compared to the rooted junction tree models in \citep{parmentier2020integer}, which is demonstrated in our computational experiments.

Against this backdrop, this paper presents how risk measures such as conditional value-at-risk and different constraints, such as chance-, budget, and logical constraints in \citep{salo2022decision} can be incorporated into the rooted junction tree model, which significantly reduces computational time compared to solving risk-averse decision strategies with decision programming. To enable our main contribution, the rooted junction tree from which the MIP is generated needs to have a specific structure. We present how rooted junction trees can be modified to achieve said structure. We also show that this can be achieved indirectly by modifying the underlying ID, which we see as more convenient from a user's standpoint.

In Section \ref{sec:background}, we present background on (LIM)IDs and the MIP reformulations of such diagrams. Section \ref{sec:contributions} continues with extending the rooted junction tree-based reformulation to consider the aforementioned risk measures and constraints. Section \ref{sec:results} presents computational results. Finally, Section \ref{sec:conclusions} concludes the paper with ideas on future research directions and the potential of reformulating IDs as MIP problems.

\section{Background}
\label{sec:background}


\subsection{Pig farm problem}
\label{subsec:pigfarm}

The pig farm problem \citep{lauritzen2001representing} is a classical example of an influence diagram and is used throughout this paper as a running example to illustrate the proposed developments. Readers should note that the pig farm problem can alternatively be cast as a partially observable Markov decision process (POMDP). \citet{cohen2023future} further discuss the modelling of POMDPs using the methodology from \citet{parmentier2020integer}. In contrast, our focus lies on the more general framework of influence diagrams.

In the pig farm problem \citep{lauritzen2001representing}, a farmer is raising pigs for a period of four months after which the pigs will be sold. During the breeding period, a pig may develop a disease, which negatively affects the retail price of the pig at the time they are sold. In the original formulation, a healthy pig commands a price of 1000 DKK and an ill pig commands a price of 300 DKK. During the first three months, a veterinarian visits the farm and tests the pigs for the disease. The specificity (or true negative rate) of the test is 80\%, whereas the sensitivity (true positive rate) is 90\%. Based on the test results, the farmer may decide to inject a medicine, which costs 100 DKK. The medicine cures an ill pig with a probability of 0.5, whereas an ill pig that is not treated is spontaneously cured with a probability of 0.1. If the medicine is given to a healthy pig, the probability of developing the disease in the subsequent month is 0.1, whereas the probability without the injection is 0.2. In the first month, a pig has the disease with a probability of 0.1.

\subsection{Influence Diagrams}
\label{subsec:ID}


An influence diagram is a directed acyclic graph $G = (N, A)$, where $N$ is the set of nodes and $A$ is the set of arcs. Let $N = N^C \cup N^D \cup N^V$ be the set of chance nodes $N^C$, value nodes $N^V$ and decision nodes $N^D$ in the ID. Let $I(j), j \in N$, denote the information set (or parents) of $j$, i.e., nodes from which there is an arc to $j$. It is typical to assume that value nodes are not parents of other nodes. IDs can intuitively represent complex decision problems with multiple periods, each containing multiple (possibly interdependent) decisions and chance nodes. For a selection of examples, see \citep{hankimaa2023decisionprogrammingjl,lauritzen2001representing,parmentier2020integer}.

Each node $j \in N$ has a discrete and finite state space $S_j$ representing possible outcomes $s_j \in S_j$. Typically, state spaces can have any discrete number of alternatives, as evidenced in the examples in \citep{hankimaa2023decisionprogrammingjl,howard2005influence,parmentier2020integer}. For a subset of nodes $C \subseteq N$, the state space and the realized outcome are defined as $S_C := \bigtimes_{c \in C}S_c$ and $s_C = (s_c)_{c \in C}$, respectively. The outcome (i.e., state) $s_j$ of a stochastic node $j \in N^C \cup N^V$ is a random variable with a probability distribution $\mathbb{P}(s_j \mid s_{I(j)})$, which corresponds to the probability of node $j$ being in state $s_j$ given that the parents $I(j)$ are in state $s_{I(j)}$. We denote $\mathbb{P}_G(s_j \mid s_{I(j)})$ when we want to emphasize that the probability distribution is associated with diagram $G$. The outcome of a decision node $d \in N^D$ is determined by a \emph{decision strategy} $\delta(s_d \mid s_{I(d)}): S_{I(d)} \bigtimes S_d \to \{0,1\}$, where $\delta(s_d \mid s_{I(d)}) = 1$ means that state $s_d$ is selected if nodes $I(d)$ are in states $s_{I(d)}$. A \emph{feasible decision strategy} is such that for each $s_{I(d)} \in S_{I(d)}$, exactly one element $s_d \in S_d$ attains $\delta(s_d \mid s_{I(d)}) = 1$ and all other $s'_d \in S_d \setminus \{s_d\}$ attain $\delta(s'_d \mid s_{I(d)}) = 0$. States $s_v \in S_v$ of a value node $v \in N^V$ represent different outcomes that have a utility value $u(s_v)$ associated with them. The total utility is calculated as a sum over the different value nodes $\sum_{v \in N^v}u(s_v)$.

The ID of the pig farm problem is presented in Figure \ref{fig:Pigs}. In the diagram, nodes $H_i$ are chance nodes representing the health status of the pig; chance nodes $T_i$ represent the test result, which is conditional on the health status of the pig; decision nodes $D_i$ represent treatment decisions; finally, value nodes $V_i$, for $i \leq 3$, represent treatment costs and the value node $V_4$ represents the pig's market price.

The solution of an ID is a decision strategy that optimizes the desired metric, typically expected utility, at the value nodes. A common additional assumption is perfect recall, meaning that previous decisions can be recalled in later stages. Under this assumption, the optimal decision strategy may be obtained by arc reversals and node removals \citep{Shachter1986} or dynamic programming \citep{TatmanShachter1990}, for example. 

Perfect recall is a rather strict assumption and in many applications, it does not hold. This challenge is circumvented with LIMIDs \citep{lauritzen2001representing}. Many algorithms for finding the decision strategy that maximizes the expected utility have been developed, such as the single policy update \citep{lauritzen2001representing}, multiple policy update \citep{MauaEtAl2012}, branch-and-bound search \citep{KhaledEtAl2013} and the aforementioned methods converting the ID to a MI(L)P \citep{parmentier2020integer,salo2022decision}.

\subsection{Rooted Junction Trees}
\label{subsec:RJT}


To achieve a MIP formulation for the decision problem, its ID, $G = (N, A)$, must first be transformed into a directed tree called a \emph{gradual rooted junction tree} (RJT) $\mathscr{G} = (\mathscr{V}, \mathscr{A})$ composed of \emph{clusters} $C \in \mathscr{V}$, which are subsets of the nodes of the ID (i.e., $C \subseteq N$) and directed arcs $\mathscr{A}$ connecting the clusters so that each cluster only has one parent. In an ID, the set of nodes $N$ consists of individual chance events, decisions and consequences, while the clusters in $\mathscr{V}$ comprise multiple nodes, hence the notational distinction between $N$ and $\mathscr{V}$. The clusters are associated with a root node $j \in N$ and we refer to clusters based on the root node as the \emph{root cluster} $C_j \in \mathscr{V}$ of node $j \in N$. Starting from an ID, we create an RJT guided by Definition \ref{def:g-rjt}, which states its necessary properties.
\begin{definition}
    \label{def:g-rjt}
    A directed rooted tree $\mathscr{G} = (\mathscr{V}, \mathscr{A})$ consisting of clusters $C_j \in \mathscr{V}$ of nodes $j \in N$ is a gradual rooted junction tree corresponding to the influence diagram $G$ if 
    \begin{enumerate}
        \item[(a)] given two clusters $C_i$ and $C_j$ in the junction tree, any cluster $C_k$ on the unique undirected path between $C_i$ and $C_j$ satisfies $C_i \cap C_j \subseteq C_k$;
        \item[(b)] each cluster $C_j \in \mathscr{V}$ is the root cluster of exactly one node $j \in N$ (that is, the root of the subgraph induced by the clusters with node $j$) and all nodes $j \in N$ appear in at least one of the clusters;
        \item[(c)] and, for each cluster, $I(j) \in C_j$.
    \end{enumerate}
\end{definition}

A rooted tree satisfying part (a) in Definition \ref{def:g-rjt} is said to satisfy the \emph{running intersection property}. As a consequence of property (a), a subgraph induced by the clusters containing node $j$ is connected. Moreover, as a consequence of property (b), we see that an RJT has as many clusters as the original influence diagram has nodes, and all nodes in the influence diagram are root nodes of exactly one cluster. Another consequence is that a cluster can contain only one node that is not contained in its parent cluster. Property (c) ensures that the root cluster contains all relevant nodes to evaluate $\delta(s_j \mid s_{I(j)})$ if $j$ is a decision node, or $\mathbb{P}(s_j \mid s_{I(j)})$ if $j$ is a chance node or a value node. 

The RJT is created by a function $f: (N,A) \to (\mathscr{V}, \mathscr{A})$. Any function that creates an RJT satisfying the properties in Definition \ref{def:g-rjt} can be used to derive the MIP formulation. In \citep{parmentier2020integer}, the authors present two alternatives for $f$. The first function uses a given topological ordering of the nodes and builds the RJT starting from the root cluster of the last node in the topological ordering and proceeding in the reverse direction of this topological ordering. This function returns an RJT with minimum treewidth given the ordering of nodes. The second function has an additional step of finding a ``good'' topological ordering that results in an RJT with minimum treewidth. For simplicity, we chose to use the function requiring a topological ordering. Using $H_1,T_1,D_1,V_1,H_2,...,H_4,V_4$ as a topological ordering, the pig farm ID in Figure \ref{fig:Pigs} is transformed to the gradual RJT in Figure \ref{fig:PigsRJT}.

\begin{figure}[ht]
\centering 
\begin{tikzpicture}
    [decision/.style={fill=white!80, draw, minimum size=2.5em, inner sep=2pt}, 
    chance/.style={circle, fill=white!80, draw, minimum size=2.5em, inner sep=2pt},
    value/.style={diamond, fill=white!60, draw, minimum size=2.5em, inner sep=2pt},
    scale=1.5]
     \node[chance] (h1) at (0, 3)      {$H_1$};
     \node[chance] (h2) at (1, 3)      {$H_2$};
     \node[chance] (h3) at (2, 3)      {$H_3$};
     \node[chance] (h4) at (3, 3)      {$H_4$};
     \node[value]  (u4) at (4, 3)      {$V_4$};
     \node[chance] (t1) at (0, 2)      {$T_{1}$};
     \node[chance] (t3) at (1, 2)      {$T_{2}$};
     \node[chance] (t5) at (2, 2)      {$T_{3}$};
     \node[decision] (d1) at (0, 1)    {$D_1$};
     \node[decision] (d2) at (1, 1)    {$D_2$};
     \node[decision] (d3) at (2, 1)    {$D_3$};
     \node[value] (u1) at (0, 0)       {$V_1$};
     \node[value] (u2) at (1, 0)       {$V_2$};
     \node[value] (u3) at (2, 0)       {$V_3$};
     \draw[->, thick] (h1) -- (t1);
     \draw[->, thick] (h1) -- (h2);
     \draw[->, thick] (h2) -- (t3);
     \draw[->, thick] (h2) -- (h3);
     \draw[->, thick] (h3) -- (t5);
     \draw[->, thick] (h3) -- (h4);
     \draw[->, thick] (h4) -- (u4);
     \draw[->, thick] (t1) -- (d1);
     \draw[->, thick] (t3) -- (d2);
     \draw[->, thick] (t5) -- (d3);
     \draw[->, thick] (d1) -- (h2);
     \draw[->, thick] (d1) -- (u1);
     \draw[->, thick] (d2) -- (h3);
     \draw[->, thick] (d2) -- (u2);
     \draw[->, thick] (d3) -- (h4);
     \draw[->, thick] (d3) -- (u3);
\end{tikzpicture}
\caption{ID of the pig farm problem \citep{lauritzen2001representing}.} \label{fig:Pigs}
\end{figure}

\begin{figure}[ht]
\centering 
\begin{tikzpicture}
    [cluster/.style={fill=white!80, draw, minimum size=2.5em, inner sep=2pt, rounded corners}]
     \node[cluster] (1) at (0, 0)      {$H_1$};
     \node[cluster] (2) at (1.5, 0)      {$H_1T_1$};
     \node[cluster] (3) at (3, 0)      {$H_1T_1D_1$};
     \node[cluster] (4) at (4.5, 0)      {$D_1V_1$};
     \node[cluster] (5) at (0, -2)      {$H_1D_1H_2$};
     \node[cluster] (6) at (1.5, -2)      {$H_2T_2$};
     \node[cluster] (7) at (3, -2)      {$H_2T_2D_2$};
     \node[cluster] (8) at (4.5, -2)      {$D_2V_2$};
     \node[cluster] (9) at (0, -4)      {$H_2D_2H_3$};
     \node[cluster] (10) at (1.5, -4)      {$H_3T_3$};
     \node[cluster] (11) at (3, -4)      {$H_3T_3D_3$};
     \node[cluster] (12) at (4.5, -4)      {$D_3V_3$};
     \node[cluster] (13) at (0, -6)      {$H_3D_3H_4$};
     \node[cluster] (14) at (1.5, -6)      {$H_4V_4$};
     \draw[->, thick] (1) -- (2);
     \draw[->, thick] (2) -- (3);
     \draw[->, thick] (3) -- (4);
     \draw[->, thick] (3) -- (5);
     \draw[->, thick] (5) -- (6);
     \draw[->, thick] (6) -- (7);
     \draw[->, thick] (7) -- (8);
     \draw[->, thick] (7) -- (9);
     \draw[->, thick] (9) -- (10);
     \draw[->, thick] (10) -- (11);
     \draw[->, thick] (11) -- (12);
     \draw[->, thick] (11) -- (13);
     \draw[->, thick] (13) -- (14);
\end{tikzpicture}
\caption{Gradual RJT of the pig farm problem.} \label{fig:PigsRJT}
\end{figure}

Formulating an optimization model based on the RJT representation starts by introducing a vector of moments $\mu_{C_j}$ for each root cluster $C_j, \ j \in N$. \citet{parmentier2020integer} show that for RJTs, we can impose constraints so that these become moments of a distribution $\mu_N$ that factorizes according to $G(N,A)$. The joint distribution $\mathbb{P}$ is said to factorize \citep{koller2009probabilistic} according to $G$ if  
\begin{equation}
    \mathbb{P}(s_N) = \prod_{j \in N}\mathbb{P}(s_j \mid s_{I(j)}). \label{eq:factorization}
\end{equation}
In the formulation, $\mu_{C_j}(s_{C_j})$ represents the probability of the nodes within the cluster $C_j$ being in states $s_{C_j}$ and part (c) of Definition \ref{def:g-rjt} ensures that $\mathbb{P}(s_j \mid s_{I(j)})$ can thus be obtained from $\mu_{C_j}(s_{C_j})$ for each $j \in N$. The resulting MIP model is 
\begin{align}
    \max &\sum_{v \in N^V} \sum_{s_{C_v} \in S_{C_v}} \mu_{C_v}(s_{C_v}) u(s_v) \label{eq:rjt-obj}\\
    \text{s.t. } & \sum_{s_{C_j} \in S_{C_j}} \mu_{C_j}(s_{C_j}) = 1, \ \forall j \in N \label{eq:rjt-probsum}\\
    & \sum_{\substack{\{s_{C_i} \in S_{C_i} \mid \\ s_{C_i \cap C_j} = s'_{C_i \cap C_j}\}}} \mu_{C_i}(s_{C_i})= \sum_{\substack{\{s_{C_j} \in S_{C_j} \mid \\ s_{C_i \cap C_j} = s'_{C_i \cap C_j}\}}} \mu_{C_j}(s_{C_j}),  \nonumber\\ &\qquad\qquad\qquad\qquad  \ \forall (C_i,C_j) \in \mathscr{A}, s'_{C_i \cap C_j} \in S_{C_i \cap C_j} \label{eq:rjt-localconsistency}\\
    & \mu_{C_j}(s_{C_j}) = \mu_{\overline{C}_j}(s_{\overline{C}_j}) \mathbb{P}(s_j \mid s_{I(j)}), \ \forall j \in N^C \cup N^V, s_{C_j} \in S_{C_j} \label{eq:rjt-prob}
    \end{align}
    \vspace{-0.7cm}
    \begin{align}
    & \mu_{C_j}(s_{C_j}) = \mu_{\overline{C}_j}(s_{\overline{C}_j})\delta(s_j \mid s_{I(j)}), \ \forall j \in N^D, s_{C_j} \in S_{C_j} \label{eq:rjt-dec}\\
    & \mu_{C_j}(s_{C_j}) \ge 0, \ \forall j \in N, s_{C_j} \in S_{C_j} \label{eq:rjt-mu}\\
    & \delta(s_j \mid s_{I(j)}) \in \{0,1\}, \ \forall j \in N^D, s_j \in S_j, s_{I(j)} \in S_{I(j)}. \label{eq:rjt-vars}
\end{align}

The formulation \eqref{eq:rjt-obj}-\eqref{eq:rjt-vars} is an expected utility maximization problem where the decision variables in the model are $\delta$ and $\mu$ and parameters are $u$ and $\mathbb{P}$. In the objective function \eqref{eq:rjt-obj}, $s_v$ is extracted from $s_{C_v}$ to evaluate the utility of each state combination of nodes in $C_v$. For notational brevity, we use $\overline{C}_j = C_j \setminus \{j\}$ to represent cluster $C_j$ without the root node $j$ in constraints \eqref{eq:rjt-prob} and \eqref{eq:rjt-dec} and $\mu_{\overline{C}_j}(s_{\overline{C}_j}) = \sum_{s_j \in S_j} \mu_{C_j}(s_{C_j})$ to represent the marginal distribution for cluster $C_j$ with the node $j$ marginalized out in constraints \eqref{eq:rjt-prob} and \eqref{eq:rjt-dec}. 

Constraints \eqref{eq:rjt-probsum} and \eqref{eq:rjt-mu} state that the variables $\mu_{C_j}$ must represent valid probability distributions, with nonnegative probabilities summing to one. Constraint \eqref{eq:rjt-localconsistency} enforces local consistency between adjacent clusters, meaning that for a pair $C_i, C_j$ of clusters connected by an edge, the marginal distribution for the nodes in both $C_i$ and $C_j$ (that is, $C_i \cap C_j$) must be the same when obtained from either $C_i$ or $C_j$. For example, for the RJT in Figure \ref{fig:PigsRJT}, constraint \eqref{eq:rjt-localconsistency} enforces that the joint probability distribution for $H_3$ and $D_3$ is the same when evaluated from clusters $H_3T_3D_3$ and $H_3D_3H_4$.

In \citep{parmentier2020integer}, the authors note that the feasibility of the decision strategy is implied by constraints \eqref{eq:rjt-probsum}-\eqref{eq:rjt-vars}. Constraint \eqref{eq:rjt-dec} enforces that $\mu_{C_j}(s_{C_j})$ either take value 0 or $\mu_{\overline{C}_j}(s_{\overline{C}_j})$, depending on the decided strategy. Constraint \eqref{eq:rjt-prob} enforces that $\mu_{C_j}(s_{C_j})$ follows the defined conditional probability distribution for $s_j$. For a more extensive explanation of the RJT model formulation, see \citep{parmentier2020integer}.

It should be noted that constraint \eqref{eq:rjt-dec} involves a product of two variables, and is thus not linear. Since we are limiting ourselves to settings with deterministic strategies (i.e., $\delta(s_d \mid s_{I(d)}): S_{I(d)} \bigtimes S_d \to \{0,1\}$), these constraints become indicator constraints and can be efficiently handled by solvers such as Gurobi \citep{gurobi}. We note that this would not be the case for more general strategies of the form $\delta(s_d \mid s_{I(d)}): S_{I(d)} \bigtimes S_d \to [0,1]$.

The number of constraints \eqref{eq:rjt-localconsistency}-\eqref{eq:rjt-dec} grows exponentially with respect to the number of nodes within a single cluster. This highlights the need to find a gradual RJT representation where the clusters are as small as possible.
\newcommand{\indep}{\perp \!\!\! \perp}
\section{Our contributions}
\label{sec:contributions}

\subsection{Extracting the utility distribution}
\label{subsec:ID-reformulation}

For problems with multiple value nodes, e.g., multi-stage decision problems, the expected utility has the property that the total expected utility is the sum of expected utilities in each value node. This property can be exploited in the solution process, and for this reason, many solution methods for IDs, including the RJT approach in \citep{parmentier2020integer}, only tackle maximum expected utility (MEU) problems. 

In contrast, risk measures (such as CVaR) require that the full probability distribution of the consequences is explicitly represented in the model. However, such representations are lost when the value nodes are placed in separate clusters, as in Figure \ref{fig:PigsRJT}, since probability distributions are only defined for each cluster separately. For example, in the pig farm problem described in Section \ref{subsec:pigfarm}, the joint distribution of $V_1$ and $V_2$ cannot be inferred from the probability distributions of clusters $C_{V_1}$ and $C_{V_2}$, as we cannot assume the probabilities of consequences in $V_1$ and $V_2$ to be independent.



The issue can be circumvented by generating the RJT based on an alternative equivalent ID $\overline{G} = (\overline{N},\overline{A})$ that collects all consequences under a single value node. We present Algorithm \ref{alg:1}, which transforms an ID into a single-value-node diagram that is equivalent to the original ID in terms of how joint probabilities and utilities are calculated. First, we formalize the notion of equivalence between IDs in Definition \ref{equivalent_diagrams}.

\begin{definition} \label{equivalent_diagrams}
We say that two IDs $G_1=(N_1,A_1)$ and $G_2=(N_2,A_2)$ are equivalent if 
\begin{enumerate}
\item[(a)] $G_1$ and $G_2$ share the same chance and decision nodes, and arcs to these nodes.
\item[(b)] There exists a bijection $g: S_{N_1} \to S_{N_2}$ such that for any $\delta(s_d \mid s_{I(d)}), \forall d \in N^{D}_1$, the following holds: $\prod_{n \in N_1}\mathbb{P}(s_n \mid s_{I(n)}) = \prod_{n \in N_2}\mathbb{P}(\bar{s}_n \mid \bar{s}_{I(n)})$ and $\sum_{v \in N_1^{v}}u(s_v) = \sum_{v \in N_2^{v}}u(\bar{s}_v)$, where $\bar{s} = g(s)$ for each $s \in S_{N_1}$.
\end{enumerate}
\end{definition}

Part (b) in Definition \ref{equivalent_diagrams} implies that for each possible state combination of value nodes from one diagram (say, $G_1$), an equivalent combination in the other diagram ($G_2$), that evaluates to the same probability and utility, must exist. For instance, if a diagram $G_1$ has two value nodes with two states each, an equivalent ID $G_2$ could have one value node with four states that correspond to all possible combinations of the two value nodes in $G_1$.

\begin{algorithm}[h]
\caption{Single-value-node conversion} \label{alg:1}
\begin{algorithmic}[1] 
\State {\bf Require} $G = (N,A)$ 
\State {\bf Initialize} $\overline{G} = (\overline{N} = \emptyset,\overline{A} =\emptyset)$ 
\State {Add $N^{C} \cup N^{D} $ to $ \overline{N}$}
\State {Set $\mathbb{P}_{\overline{G}}(s_{j} \mid s_{I(j)}) = \mathbb{P}_{G}(s_{j} \mid s_{I(j)}), \forall j \in N^{C}$}
\State {Add $\{(a,b) \in A \mid b \in N^{C} \cup N^{D}\} $ to $\overline{A}$} 
\State {Create node $\overline{v}$ such that $S_{\overline{v}} = \bigtimes_{v \in N^{v}}S_v$, add $\overline{v} \in \overline{N}$}
\State {Add $(a,\overline{v}) $ to $ A$ for each $a \in I(v)$ such that $v \in N^{V}$}
\State {Set $u(s_{\overline{v}}) = \sum_{v \in N^V}u(s_{v})$}
\State {Set $\mathbb{P}( s_{\overline{v}} \mid s_{I(\overline{v})}) = \prod_{v \in N^{v}}\mathbb{P}(s_{v} \mid s_{I(v)})$}
\State{\bf Return} $\overline{G}$
\end{algorithmic}
\end{algorithm}

\begin{proposition}
Let $G = (N,A)$ be an ID with $|N^{V}| > 1$. An ID $\overline{G} = (\overline{N},\overline{A})$ constructed using Algorithm \ref{alg:1} is equivalent (cf. Definition \ref{equivalent_diagrams}) . 
\end{proposition}

\begin{proof}
See \ref{sec:appendix_in_text}.
\end{proof}

When applying Algorithm \ref{alg:1} and transforming the ID into an RJT, according to Definition \ref{def:g-rjt}, part (c), we have that $\bigcup_{v \in V}I(v) \subseteq C_{\bar{v}}$, where $\bar{v}$ represents the unique value node in the new diagram. Consequently, the marginal probability distribution $\mu_{C_{\bar{v}}}$ contains information on the joint probability distribution of the consequences $\mathbb{P}( s_{\overline{v}} \mid s_{I(\overline{v})})$ and this can be used to expose the probability distribution of the utility values. Following this approach, the modified ID of the pig farm problem is presented in Figure \ref{fig:PigsReformulated} and the corresponding gradual RJT in Figure \ref{fig:PigsRJTReformulated}.

However, this incurs in computationally more demanding versions of model \eqref{eq:rjt-obj}-\eqref{eq:rjt-vars}. In the single-value-node version of the pig farm problem, all decision nodes $D_k$, $k=1,2,3$, are in the information set of $\bar{V}$. It follows from the running intersection property that $D_k$ must be contained in every cluster that is in the undirected path between $C_{D_k}$ and $C_{\bar{v}}$. Therefore, the clusters become larger as the parents of value nodes are ``carried over'', instead of evaluating separable components of the utility function at different value nodes. As discussed in \citep{parmentier2020integer}, this increases the computational complexity of the resulting model.


\begin{figure}[ht]
\centering 
\begin{tikzpicture}
    [decision/.style={fill=white!80, draw, minimum size=2.5em, inner sep=2pt}, 
    chance/.style={circle, fill=white!80, draw, minimum size=2.5em, inner sep=2pt},
    value/.style={diamond, fill=white!60, draw, minimum size=2.5em, inner sep=2pt},
    scale=1.5]
     \node[chance] (h1) at (0, 4)      {$H_1$};
     \node[chance] (h2) at (1, 4)      {$H_2$};
     \node[chance] (h3) at (2, 4)      {$H_3$};
     \node[chance] (h4) at (3, 4)      {$H_4$};
     \node[value]  (u4) at (4, 4)      {$\bar{V}$};
     \node[chance] (t1) at (0, 3)      {$T_{1}$};
     \node[chance] (t3) at (1, 3)      {$T_{2}$};
     \node[chance] (t5) at (2, 3)      {$T_{3}$};
     \node[decision] (d1) at (0, 1.4)    {$D_1$};
     \node[decision] (d2) at (1, 1.4)    {$D_2$};
     \node[decision] (d3) at (2, 1.4)    {$D_3$};
     \draw[->, thick] (h1) -- (t1);
     \draw[->, thick] (h1) -- (h2);
     \draw[->, thick] (h2) -- (t3);
     \draw[->, thick] (h2) -- (h3);
     \draw[->, thick] (h3) -- (t5);
     \draw[->, thick] (h3) -- (h4);
     \draw[->, thick] (h4) -- (u4);
     \draw[->, thick] (t1) -- (d1);
     \draw[->, thick] (t3) -- (d2);
     \draw[->, thick] (t5) -- (d3);
     \draw[->, thick] (d1) -- (h2);
     \draw[->, thick] (d2) -- (h3);
     \draw[->, thick] (d3) -- (h4);
     \draw[->, thick] (d1.east) -- (u4);
     \draw[->, thick] (d2.east) -- (u4);
     \draw[->, thick] (d3.east) -- (u4);
\end{tikzpicture}
\caption{The pig farm problem reformulated ID.} \label{fig:PigsReformulated}
\end{figure}

\begin{figure}[ht]
\centering 
\begin{tikzpicture}
    [cluster/.style={fill=white!80, draw, minimum size=2.5em, inner sep=2pt, rounded corners}]
     \node[cluster] (1) at (0, 0)      {$H_1$};
     \node[cluster] (2) at (1.5, 0)      {$H_1T_1$};
     \node[cluster] (3) at (3, 0)      {$H_1T_1D_1$};
     \node[cluster] (5) at (0, -1.5)      {$H_1D_1H_2$};
     \node[cluster] (6) at (2, -1.5)      {$D_1H_2T_2$};
     \node[cluster] (7) at (4, -1.5)      {$D_1H_2T_2D_2$};
     \node[cluster] (9) at (0, -3)      {$D_1H_2D_2H_3$};
     \node[cluster] (10) at (2.5, -3)      {$D_1D_2H_3T_3$};
     \node[cluster] (11) at (5, -3)      {$D_1D_2H_3T_3D_3$};
     \node[cluster] (13) at (0, -4.5)      {$D_1D_2H_3D_3H_4$};
     \node[cluster] (14) at (3, -4.5)      {$D_1D_2D_3H_4\bar{V}$};
     \draw[->, thick] (1) -- (2);
     \draw[->, thick] (2) -- (3);
     \draw[->, thick] (3) -- (5);
     \draw[->, thick] (5) -- (6);
     \draw[->, thick] (6) -- (7);
     \draw[->, thick] (7) -- (9);
     \draw[->, thick] (9) -- (10);
     \draw[->, thick] (10) -- (11);
     \draw[->, thick] (11) -- (13);
     \draw[->, thick] (13) -- (14);
\end{tikzpicture}
\caption{Gradual RJT of the reformulated pig farm problem.} \label{fig:PigsRJTReformulated}
\end{figure}


\subsection{Modifying the RJT}
\label{sec:mod_rjt}

A single-value-node ID guarantees that the full probability distribution of the consequences can be exposed in the MIP model to optimize or constrain different risk metrics. However, decision-makers may not only be interested in imposing risk constraints on consequences. For instance, in the pig farm problem, the farmer may want to impose chance constraints to ensure that pigs stay healthy throughout the breeding period with a certain likelihood. This requires that the joint probability distribution of the health states of a pig in all periods is available. That is, one needs an RJT that has a cluster containing all health status nodes $H_1,...,H_4$, which is not found in the RJTs created with functions from \citep{parmentier2020integer} or with Algorithm \ref{alg:1} (Figures \ref{fig:PigsRJT} and \ref{fig:PigsRJTReformulated}).

A way to generate a cluster that exposes the desired probability distribution is to directly modify an RJT obtained following \citep{parmentier2020integer} while ensuring that the modified RJT fulfills Definition \ref{def:g-rjt}. Assume that there exists a topological ordering among the nodes of the RJT $\mathscr{G} = (\mathscr{V}, \mathscr{A})$. Suppose one wants to create a cluster that exposes the joint probability distribution of nodes $M \subseteq N$. Then, one needs to generate an RJT $\bar{\mathscr{G}} = (\bar{\mathscr{V}}, \bar{\mathscr{A}})$ such that $\exists C'_v \in \bar{\mathscr{V}} $ such that $ M \subseteq C'_v$. Let $\preccurlyeq$ represent a topological order of the nodes $N$ and $\max_{\preccurlyeq}M$ return the node with the highest topological order in set $M$. Let \begin{eqnarray*}
    P(C_1,C_k) := \{C_j \in \mathscr{V} \mid \exists (C_1,...,C_j,...,C_k) \\ \text{ with } (C_l,C_{l+1}) \in \mathscr{A}, \forall l \in \{1,\dots,k\}\}
\end{eqnarray*} be the set of clusters that are contained on any directed path between clusters $C_1$ and $C_k$ (including themselves). Let $F(j) := \{k \in N \mid P(C_j,C_k) \neq \emptyset\}$ be the set of nodes, whose root cluster can be reached via a directed path starting from cluster $C_j$. Then, Algorithm \ref{alg:2} leads to an RJT with the desired structure.

\begin{algorithm}[h]
\caption{RJT modification} \label{alg:2}
\begin{algorithmic}[1]
\State {\bf Require}  $M \subseteq N$, Topological order $\preccurlyeq$
\State {\bf Initialize} $\bar{\mathscr{G}} = (\bar{\mathscr{V}}, \bar{\mathscr{A}})$ equal to $\mathscr{G}$ and denote $C'_n$ as the root cluster of $n \in N$ in $\bar{\mathscr{G}}$
\State {Find $m = \max_\preccurlyeq M$ }
\State {{\bf For } $n \in M\setminus\{m\}$ such that $n \notin C_m$ {\bf do:}}
\Indent
\State{{\bf If }$m \notin F(n)$}
\Indent
\State {Find $e = max_{\preccurlyeq} \{j \in N \mid n,m \in F(j)\}$}
\State {Find $g \in N \text{ such that } (C_e,C_g) \in \mathscr{A}, m \in F(g), \linebreak n \notin F(g)$}
\State {Set $C'_e \cap C'_g \in C'_c, \forall C'_c \in P(C_e,C_n)$}
\State {Set $(C'_e,C'_g) \notin \bar{\mathscr{A}}$, $(C'_n,C'_g) \in \bar{\mathscr{A}}$}
\EndIndent
\State {Set $n \in C'_a, \forall C'_a \in P(C'_n,C'_m)$}
\EndIndent
\State {\bf Return $\bar{\mathscr{G}}$}
\end{algorithmic}
\end{algorithm}

\begin{proposition}
Assume that an RJT $\mathscr{G} = (\mathscr{V}, \mathscr{A})$ satisfies Definition \ref{def:g-rjt}. The RJT $\bar{\mathscr{G}} = (\bar{\mathscr{V}}, \bar{\mathscr{A}})$ generated from $\mathscr{G}$ using Algorithm \ref{alg:2} satisfies Definition \ref{def:g-rjt}.
\end{proposition}

\begin{proof}
See \ref{sec:appendix_in_text}.
\end{proof}

As an example, we can apply Algorithm \ref{alg:2}  to the RJT in Figure \ref{fig:PigsRJT}, which is created by a function given in \citep{parmentier2020integer}, to include a cluster that contains nodes $H_1,...,H_4$. The output of Algorithm \ref{alg:2} is the diagram presented in Figure \ref{fig:PigsRJT_mod}. The cluster that exposes the desired joint probability distribution for nodes $H_1,...,H_4$ is the root cluster of $H_4$. A more detailed example of applying Algorithm \ref{alg:2} can be found in \ref{sec:appendix_in_text}.

\begin{figure}[ht]
\centering 
\begin{tikzpicture}
    [cluster/.style={fill=white!80, draw, minimum size=2.5em, inner sep=2pt, rounded corners}]
     \node[cluster] (1) at (0, 0)      {$H_1$};
     \node[cluster] (2) at (1.5, 0)      {$H_1T_1$};
     \node[cluster] (3) at (3, 0)      {$H_1T_1D_1$};
     \node[cluster] (4) at (5, 0)      {$D_1V_1$};
     \node[cluster] (5) at (0, -2)      {$H_1D_1H_2$};
     \node[cluster] (6) at (2, -2)      {$H_1H_2T_2$};
     \node[cluster] (7) at (4, -2)      {$H_1H_2T_2D_2$};
     \node[cluster] (8) at (6, -2)      {$D_2V_2$};
     \node[cluster] (9) at (0, -4)      {$H_1H_2D_2H_3$};
     \node[cluster] (10) at (2.5, -4)      {$H_1H_2H_3T_3$};
     \node[cluster] (11) at (5, -4)      {$H_1H_2H_3T_3D_3$};
     \node[cluster] (12) at (7, -4)      {$D_3V_3$};
     \node[cluster] (13) at (0, -6)      {$H_1H_2H_3D_3H_4$};
     \node[cluster] (14) at (2, -6)      {$H_4V_4$};
     \draw[->, thick] (1) -- (2);
     \draw[->, thick] (2) -- (3);
     \draw[->, thick] (3) -- (4);
     \draw[->, thick] (3) -- (5);
     \draw[->, thick] (5) -- (6);
     \draw[->, thick] (6) -- (7);
     \draw[->, thick] (7) -- (8);
     \draw[->, thick] (7) -- (9);
     \draw[->, thick] (9) -- (10);
     \draw[->, thick] (10) -- (11);
     \draw[->, thick] (11) -- (12);
     \draw[->, thick] (11) -- (13);
     \draw[->, thick] (13) -- (14);
\end{tikzpicture}
\caption{Gradual RJT for pig farm problem with a cluster that contains nodes $H_1,...,H_4$.} \label{fig:PigsRJT_mod}
\end{figure}

\subsection{Imposing chance, logical, and budget constraints}
\label{sect:chance_and_budget_cn}

Our proposed developments allow one to expose the joint probability distribution of any combination of nodes in the ID, which in turn, enables the formulation of a broad range of risk-aversion-related constraints. 

Chance constraints for the joint probability of nodes $M \subseteq N$ can be imposed on the marginal probability distribution of a cluster $C_n$ such that $M \subseteq C_n$. If no such cluster exists in the generated RJT, a suitable cluster can be created with Algorithm \ref{alg:2}, or regenerating the RJT based on an ID created with Algorithm \ref{alg:1} if $M$ only contains value nodes and their parents. Then, chance constraints can be imposed as follows:
\begin{equation}
    \sum_{s_{C_{n}} \in S_{C_{n}} \mid s_{n} \in S^o_{{n}}}\mu(s_{C_{n}}) \leq p, \label{eq:chance_constraint}
\end{equation}
where $S^o_{{n}}$ is the set of outcomes that the decision maker wishes to constrain and $p \in [0,1]$ represents a threshold. For instance, assume that a decision-maker wishes to enforce that the probability of the payout of the process being less than some fixed limit $b$ is at most $p$. Then, a suitable RJT can be generated based on an ID created with Algorithm \ref{alg:1}. Constraints can then be enforced for the root cluster of the single value node $C_{\bar{v}}$ and $S^o_{\bar{v}}$ would contain all states $s_{\bar{v}}$ such that $u(s_{\bar{v}}) < b$. Note that this formulation can be enforced for any cluster $C_k$ such that $M \subseteq C_k$.

Logical constraints can be seen as a special case of chance constraints. For example, in the pig farm problem (in Section \ref{sec:background}), the farmer may wish to attain an optimal decision strategy while ensuring that the number of injections is at most two per pig due to, e.g., potential side effects. Then, $S^{o}_{\bar{v}}$ would contain all realizations of the nodes in $C_{\bar{v}}$ that would lead to a violation of the constraint, i.e., the state combinations in which three injections would be given to a pig. In that case, constraint \eqref{eq:logical_constraint} that makes these scenarios infeasible could be imposed. 
\begin{equation}
    \sum_{s_{C_{\bar{v}}} \in S_{C_{\bar{v}}} \mid s_{{\bar{v}}} \in S^o_{{\bar{v}}}}\mu(s_{C_{\overline{v}}}) \leq 0. \label{eq:logical_constraint}
\end{equation}

Budget constraints are analogous to logical constraints, as the farmer could instead have an injection budget, say 200 DKK per pig. Then, $S^{o}_{\bar{v}}$ should contain all states $s_{\bar{v}}$, where more than 200 DKK is used for treating a pig, with the constraint enforced similarly as in (\ref{eq:logical_constraint}).
\subsection{Conditional Value-at-Risk}\label{sect:Cvar}
In addition to a number of risk constraints, the proposed reformulations also enable the consideration of alternative risk measures. Next, we focus our presentation on how to maximize conditional value-at-risk. However, we highlight that other risk metrics such as absolute or lower semi-absolute deviation \citep{salo2022decision} can, in principle, be used. The entropic risk measure \citep{follmer2011entropic} can also be used as a constraint. However, incorporating it in the objective function will introduce nonlinearity in the model due to the logarithmic nature of the measure.

The proposed formulation for conditional value-at-risk maximization is analogous to the method developed for decision programming in \citep{salo2022decision}. It assumes that the joint probability distribution of utility values is available, and hence an RJT generated based on the single-value-node representation of the ID is sufficient for generating a suitable cluster. Let us assume that the decision problem has a single value node $\bar{v}$ with possible utility values $u \in U$. Let $p(u)$ be the probability of attaining utility value $u$. In the presence of a single value node, we would define $p(u) = \sum_{s_{C_{\bar{v}}} \in S_{C_{\bar{v}}} | U(s_{C_{\bar{v}}}) = u} \mu(s_{C_{\bar{v}}})$ and pose the constraints 
\begin{align}
    &\eta - u \le M \lambda(u), \ & \forall u \in U \label{eq:cvar-1}\\ 
    &\eta - u \ge (M+\epsilon) \lambda(u) - M, \ & \forall u \in U \label{eq:cvar-2}\\ 
    &\eta - u \le (M+\epsilon)\overline{\lambda}(u) - \epsilon, \ & \hspace{1cm}\forall u \in U \label{eq:cvar-3}\\ 
    &\eta - u \ge M (\overline{\lambda}(u)-1), \ & \forall u \in U \label{eq:cvar-4}\\ 
    &\overline{\rho}(u) \le \overline{\lambda}(u), \ & \forall u \in U \label{eq:cvar-5}
    \end{align}
    \vspace{-0.95cm}
    \begin{align}
    p(u) - (1-\lambda(u)) &\le \rho(u) \le \lambda(u) , \ & \forall u \in U \label{eq:cvar-6}
    \end{align}
    \vspace{-0.95cm}
    \begin{align}
    &\rho(u) \le \overline{\rho}(u) \le p(u), \ & \hspace{2cm}\forall u \in U \label{eq:cvar-7}\\ 
    &\sum_{u \in U} \overline{\rho}(u) = \alpha \label{eq:cvar-8}
    \end{align}
    \vspace{-0.8cm}
    \begin{align}
    &\lambda(u), \overline{\lambda}(u) \in \{0,1\}, \ & \hspace{2.1cm}\forall u \in U \label{eq:cvar-9}\\ 
    &\rho(u), \overline{\rho}(u) \in [0,1], \ & \forall u \in U \label{eq:cvar-10}\\
    &\eta \in \mathbb{R}, \label{eq:cvar-11}
\end{align}
where $\alpha$ is the probability level in VaR$_{\alpha}$.  Table \ref{tbl:cvar-vars} describes which values the decision variables take due to the constraints \eqref{eq:cvar-1}-\eqref{eq:cvar-11}.
\begin{table}[ht]
\centering
\begin{tabular}{l|l}
Variable            & Value                                 \\ \hline
$\eta$              & VaR$_{\alpha}$                        \\
$\lambda(u)$        & 1 if $u < \eta$                       \\
$\overline{\lambda}(u)$  & 0 if $u > \eta$                       \\
$\rho(u)$           & 0 if $\lambda(u)=0$, $p(u)$ otherwise \\
$\overline{\rho}(u)$     & $\begin{cases}
    p(u) &\text{ if } u<\eta,\\
    \alpha-\sum_{u \in U} p(u) &\text{ if } u=\eta, \\
    0 &\text{ if } u>\eta \ (\bar{\lambda}(u)=0)
\end{cases} $
\end{tabular}
\caption{Variables and the corresponding values that satisfy \eqref{eq:cvar-1}-\eqref{eq:cvar-10}}
\label{tbl:cvar-vars}
\end{table}

In constraints \eqref{eq:cvar-1}-\eqref{eq:cvar-10}, $M$ is a large positive number and $\epsilon$ is a small positive number. The parameter $\epsilon$ is used to model strict inequalities, which cannot be directly used in mathematical optimization solvers. For example, $x \ge \epsilon$ is assumed to be equivalent to $x > 0$. In practice, it is enough to set $\epsilon$ strictly smaller than the minimum difference of distinct utility values. In \citep{salo2022decision}, the authors use $\epsilon = \frac{1}{2}\min\{|U(s_{\overline{v}}) - U(s'_{\overline{v}})| : |U(s_{\overline{v}}) - U(s'_{\overline{v}})| > 0, s_{\overline{v}}, s'_{\overline{v}} \in S_{\overline{v}}\}$. When $\lambda(u) = 0$, constraints \eqref{eq:cvar-1} and \eqref{eq:cvar-2} become $-M \le \eta - u \le 0$, or $\eta \le u$. When $\lambda(u) = 1$, they instead become $\epsilon \le \eta - u \le M$, or $\eta > u$. Constraints \eqref{eq:cvar-3} and \eqref{eq:cvar-4} can be examined similarly to obtain the results in Table \ref{tbl:cvar-vars}.

The correct behavior of variables $\rho(u)$ is enforced by \eqref{eq:cvar-6} and \eqref{eq:cvar-7}. If $\lambda(u)=0$, constraint \eqref{eq:cvar-6} forces $\rho(u)$ to zero. If $\lambda(u)=1$, then $\rho(u)=p(u)$. Finally, assuming $\eta$ is equal to VaR$_{\alpha}$ and $\rho(u)$ equal to $p(u)$ for all $u < \eta$, the value of $\overline{\rho}(u)$ must be $\alpha-\sum_{u \in U} \rho(u)$ for $u=\eta$. It is easy to see that $\eta$ must be equal to VaR$_{\alpha}$ for there to be a feasible solution for the other variables. For an equivalence proof, see \citet[][Appendix A]{salo2022decision}.

By introducing constraints \eqref{eq:cvar-1}-\eqref{eq:cvar-10}, the $CVaR_{\alpha}$ can then be obtained as $\frac{1}{\alpha}\sum_{u \in U}\overline{\rho}(u)u$. This can be either used as in the objective function or as a part of the constraints of the problem. We also note that the described approach is very versatile in that $u$ can be selected to be, e.g., a stage-specific utility function, thus allowing us to limit risk in specific stages of a multi-stage problem. \citet{krokhmal2002portfolio} discusses the implications of stage-wise CVaR-constraints in detail.
\section{Computational experiments}
\label{sec:results}

To assess the computational performance of the model \eqref{eq:rjt-obj}-\eqref{eq:rjt-vars}, we use the pig farm problem described earlier. We compare two different versions of the pig farm problem: one with the RJT formulation and the other with the decision programming formulation from \citep{hankimaa2023decisionprogrammingjl}. An additional computational example and an analysis of the resulting model sizes from each formulation can be found in \ref{sec:appendix_in_text}. All problems were solved using a single thread on an Intel E5-2680 CPU at 2.5GHz and 16GB of RAM, provided by the Aalto University School of Science "Science-IT" project. The models were implemented using Julia v1.10.3 \citep{Julia-2017} and JuMP v1.23.0 \citep{DunningHuchetteLubin2017} and solved with the Gurobi solver v11.0.2 \citep{gurobi}. The code and data used in this section are available at \href{https://github.com/gamma-opt/risk-averse-RJT}{www.github.com/gamma-opt/risk-averse-RJT}.

\subsection{Risk-averse pig farm problem}

A risk-averse version of the pig farm problem, which maximizes conditional value-at-risk for a confidence level of $1-\alpha=0.85$ is solved for different numbers of breeding periods. We use the RJT based on the single-value-node ID from Figure \ref{fig:PigsRJTReformulated} to create the optimization model and add the constraints described in Section \ref{sect:Cvar} to represent CVaR. The solution times using our RJT-based formulation are compared to the solution times derived using the decision programming formulation from \citep{hankimaa2023decisionprogrammingjl}. For practical reasons, we solve the same single-value-node version of the pig farm problem to compare the solutions. In practice, decision programming can maximize CVaR in IDs with any number of value nodes. However, decision programming creates the exact same MILP model regardless of the number of value nodes.  

\begin{figure}[ht]
    \centering
    \includegraphics[width=0.45\textwidth]{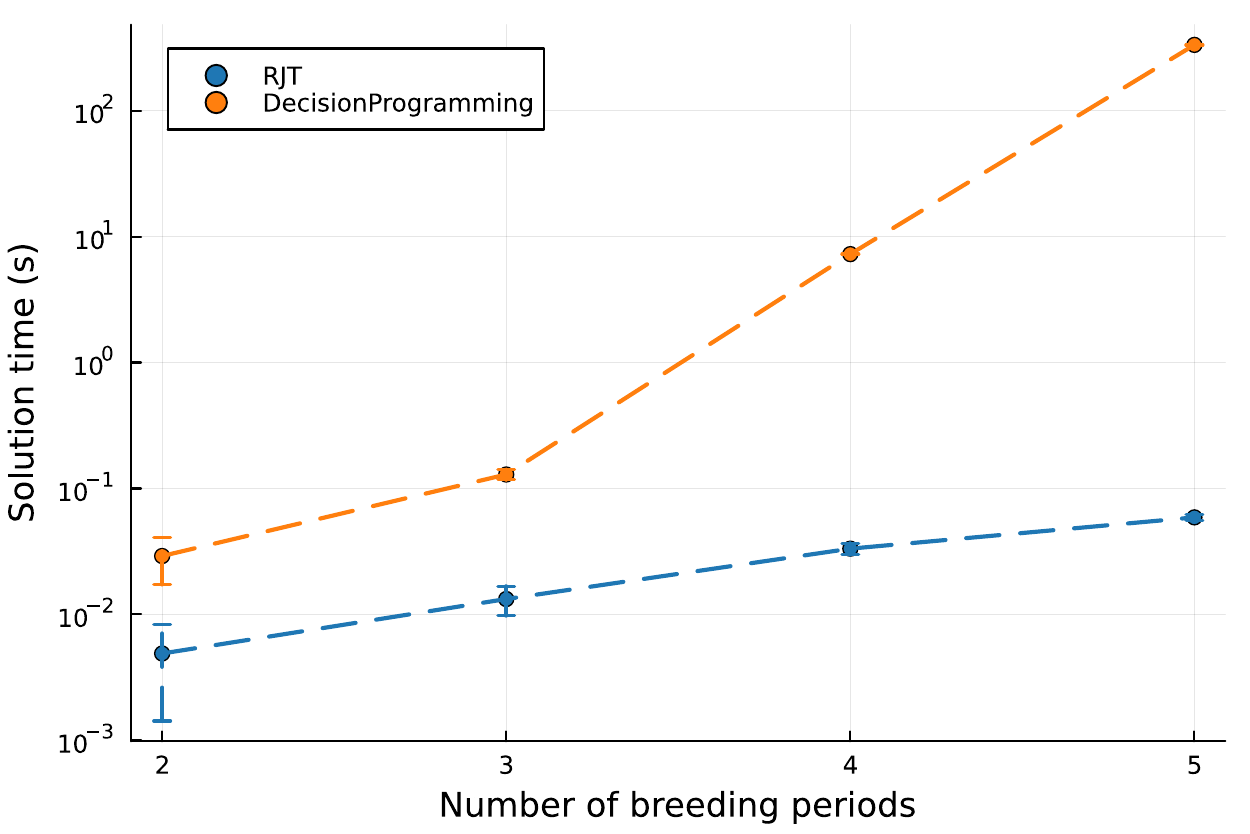}
    \caption{Mean solution times and standard deviations (bars) for 50 random instances in the risk-averse pig farm problem with 2-5 breeding periods on a logarithmic scale.}
    \label{fig:sol_times_DP_vsRJT}
\end{figure}

The solution times of 50 randomly generated instances of the risk-averse pig farm problem with different sizes are presented in Figure \ref{fig:sol_times_DP_vsRJT}. The RJT-based formulation consistently offers better computational performance than decision programming for the pig farm problem. In larger pig farm instances, the RJT-based formulation is three orders of magnitude faster than the decision programming formulation. However, the RJT-based formulation still grows exponentially with respect to the number of breeding periods, which could result in computational challenges for larger instances. Still, this exponential growth of the RJT model can be seen as a worst-case scenario, while many problems, including the original pig farm problem (Figure \ref{fig:Pigs}), exhibit treewidth independent of the number of stages. For decision programming, as discussed in \citep{hankimaa2023decisionprogrammingjl}, the number of constraints is exponential in the number of nodes.

\subsection{Chance-constrained pig farm problem}

\begin{figure}[ht]
    \centering
    \includegraphics[width=0.45\textwidth]{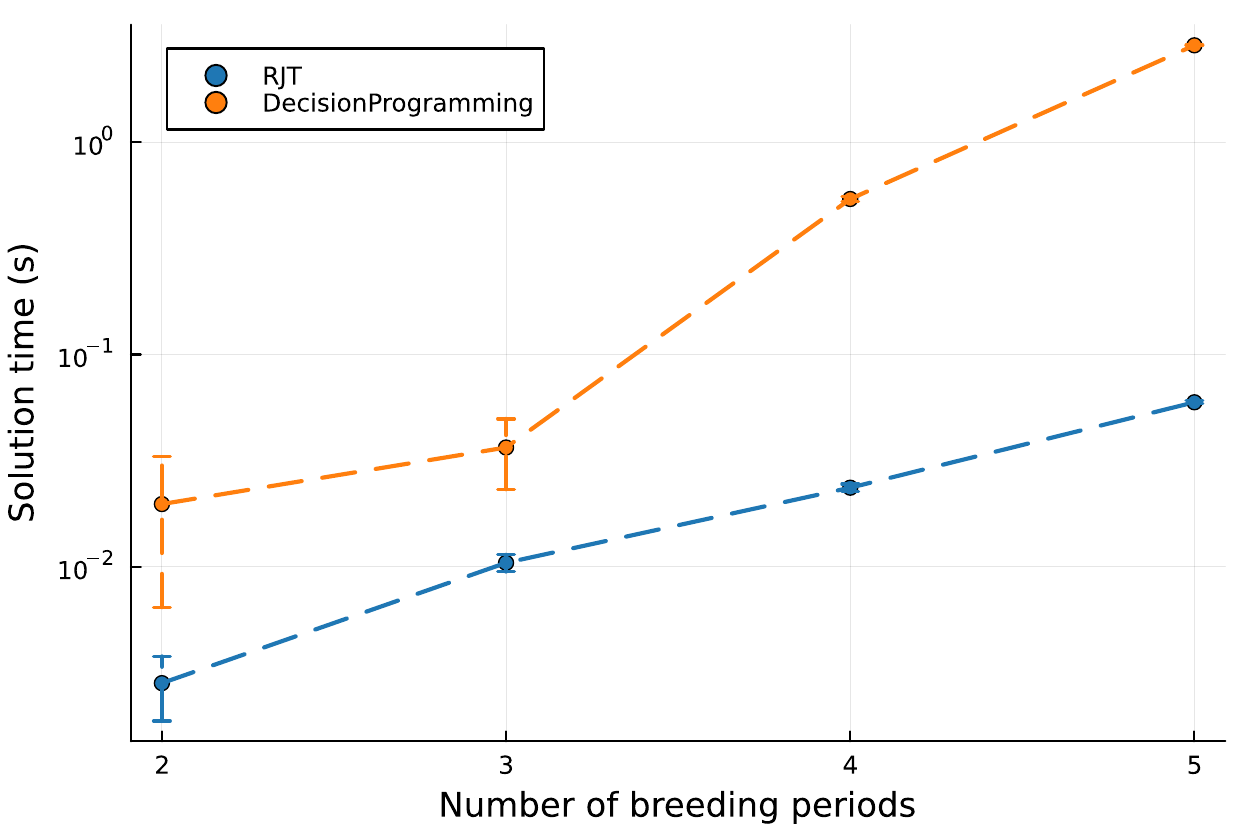}
    \caption{Mean solution times and standard deviations (bars) for 50 random instances in the chance-constrained pig farm problem with 2-5 breeding periods on a logarithmic scale.}
    \label{fig:sol_times_DP_vsRJTCC}
\end{figure}

In addition, we analyze the computational performance of our formulations on a chance-constrained version of the pig farm problem with different numbers of breeding periods. We use the RJT in Figure \ref{fig:PigsRJT_mod} and assign chance constraints to the root cluster of $H_4$ enforcing that the probability of a pig being ill at any time during the breeding period must be less than $40~\%$. Chance constraints are enforced as described in Section \ref{sect:chance_and_budget_cn}. In Figure \ref{fig:sol_times_DP_vsRJTCC}, we compare the results by solving the same problem with decision programming \citep{hankimaa2023decisionprogrammingjl}.

The optimization model created based on RJT solves the chance-constrained problem faster than the corresponding decision programming model. In accordance with the results of the risk-averse pig farm problem, RJT is an order of magnitude faster than the corresponding decision programming model.
\section{Conclusions}
\label{sec:conclusions}
\vspace{-0.1cm}

In this paper, we have described a MIP reformulation of decision problems presented as IDs, originally proposed in \citep{parmentier2020integer}. Our main contribution is to extend the modelling framework proposed by \citet{parmentier2020integer} to embed it with more general modelling capabilities. We illustrate how chance constraints and conditional value-at-risk can be incorporated into the formulation. We demonstrate how suitable RJTs can be generated, either by modifying the underlying ID (Algorithm \ref{alg:1}) or directly modifying the RJT (Algorithm \ref{alg:2}).

We show that the model in \citep{parmentier2020integer} can be extended beyond expected utility maximization problems to incorporate most of the constraints and objective functions present in decision programming, the alternative MIP reformulation based on LIMIDs described by \citet{salo2022decision} and \citet{hankimaa2023decisionprogrammingjl}. The advantage of using the models described in this paper is that in terms of model size, decision programming models grow exponentially with respect to the number of nodes, whereas the RJT model grows exponentially with respect to treewidth, which is only indirectly influenced by the number of nodes.

We also present computational results comparing the computational performance of decision programming and our extension of the RJT model when applied to risk-averse and chance-constrained variants of the pig farm problems. The computational results indicate that risk-averse decision strategies for IDs can be solved considerably faster by using the RJT formulation.

Although this paper furthers the state-of-the-art for MIPs solving risk-averse IDs, typically the resulting MIP is a large-scale model, which in turn limits the size of the problems that can be solved. Hence, future research should concentrate on improving the computational tractability of the MIP model by developing specialized decomposition methods and more efficient formulations.
\appendix
\appendix
\section{Supplementary material}
\label{sec:appendix_in_text}

\section*{Acknowledgements}
This work was supported by the Research Council of Finland (decision 332180). We are also thankful for the contributions from Prof. Ahti Salo and the computer resources from the Aalto University School of Science “Science-IT” project.




\bibliographystyle{elsarticle-harv} 
\bibliography{RJT_bibliography}


\end{document}


\begin{frontmatter}
\title{Supplementary material for ``Risk-averse decision strategies for influence diagrams using rooted junction trees'' }

\author[label1]{Olli Herrala}
\author[label1]{Topias Terho}
\author[label1]{Fabricio Oliveira\corref{cor1}}
\cortext[cor1]{Corresponding author: fabricio.oliveira@aalto.fi}
\affiliation[label1]{organization={Department of Mathematics and Systems Analysis},
            addressline={Aalto University, School of Science}, 
            postcode={FI-00076 Aalto}, 
            country={Finland}}
\end{frontmatter}

\section{Introduction}

This supplementary material contains a more detailed description of Algorithm 2, a discussion on the size of model (2)-(8), the proofs for Propositions 3.2 and 3.3, as well as an additional computational example.

\section{Example of applying Algorithm 2}

Algorithm 2 is explained using the example RJT presented in Figure \ref{fig:rjt_example}. Assume that the topological order of the nodes is given as follows: $\{A,B,C,D,E,F\}$ and assume that the aim is to modify the RJT by creating a cluster containing nodes $A, E$ and $F$.

First, Algorithm 2 initializes a new RJT containing the same nodes and arcs. Line 2 of Algorithm 2 finds the last node in $M$ according to the topological order. In the example diagram, it is node $F$. After running Algorithm 2, cluster $C_F$, which currently contains nodes $B$ and $F$, will be the cluster containing nodes $A, E$ and $F$. Line 3 of the algorithm checks if any of the nodes are already contained in $C_F$. In the example RJT, out of $A, E$ and $F$, only $F$ is contained in $C_F$. The algorithm will then run the for-loop starting in Line 4 for both $A$ and $E$. First, consider $n = A$. In Line 5, Algorithm 2 checks if there exists a directed path between $C_A$ and $C_F$. In the example diagram, such path  $p = (C_A,C_B,C_D,C_F)$ exists. Therefore, Algorithm 2 proceeds directly to Line 10, in which $A$ is added to all clusters in path $p$. This results in the RJT in Figure \ref{fig:rjt_example2}. 

Next, the algorithm switches to $n = E$. Now there is no directed path from $C_E$ to $C_F$ in the RJT and therefore Algorithm 2 conducts Lines 6-9. In Line 6, Algorithm 2 searches for a node and its root cluster from which $C_F$ and $C_E$ both can be reached. In the example, there are two potential nodes, $A$ and $B$. Then, the Algorithm selects the node with maximum topological order, which in this case is $B$. In Line 7, Algorithm 2 finds a direct predecessor of $C_B$, which is in a path towards $C_F$. In the example, the resulting node is $D$ and the corresponding cluster is $C_D$. This step returns a unique solution because $F$ can only be in one such branch as it has to have a unique root cluster. In Line 8, the algorithm fills clusters on the path between $C_B$ and $C_E$ with nodes contained in both clusters $C_B$ and $C_D$. This corresponds to nodes $A$ and $B$. Both are already contained in cluster $C_C$, whereas $A$ must be added to $C_E$. This creates an RJT in Figure \ref{fig:rjt_example3}. In Line 9, Algorithm 2 removes the arc between clusters $C_B$ and $C_D$ and adds an arc between clusters $C_E$ and $C_D$ to create the RJT in Figure \ref{fig:rjt_example4}. The procedure in Lines 6-9 has now created an RJT with a directed path between clusters $C_E$ and $C_F$. Finally, Algorithm 2 implements Line 10 and adds $E$ into $C_D$ and $C_F$ to create the RJT in Figure \ref{fig:rjt_example5}. Then, Algorithm 2 returns an RJT where cluster $C_F$ contains the desired nodes.

\begin{figure}[ht]
\centering 
\begin{tikzpicture}
    [cluster/.style={fill=white!80, draw, minimum size=2.5em, inner sep=2pt, rounded corners}]
     \node[cluster] (1) at (0, 2)      {$C_A: A$};
     \node[cluster] (2) at (2, 2)      {$C_B: AB$};
     \node[cluster] (3) at (4, 1)      {$C_C: ABC$};
     \node[cluster] (5) at (6, 1)      {$C_E: BCE$};
     \node[cluster] (6) at (4, 3)      {$C_D: BD$};
     \node[cluster] (7) at (6, 3)      {$C_F: BF$};
     \draw[->, thick] (1) -- (2);
     \draw[->, thick] (2) -- (3);
     \draw[->, thick] (2) -- (6);
     \draw[->, thick] (3) -- (5);
     \draw[->, thick] (6) -- (7);

\end{tikzpicture}
\caption{An example RJT.} \label{fig:rjt_example}
\end{figure}

\begin{figure}[ht]
\centering 
\begin{tikzpicture}
    [cluster/.style={fill=white!80, draw, minimum size=2.5em, inner sep=2pt, rounded corners}]
     \node[cluster] (1) at (0, 2)      {$C_A: A$};
     \node[cluster] (2) at (2, 2)      {$C_B: AB$};
     \node[cluster] (3) at (4, 1)      {$C_C: ABC$};
     \node[cluster] (5) at (6, 1)      {$C_E: BCE$};
     \node[cluster] (6) at (4, 3)      {$C_D: ABD$};
     \node[cluster] (7) at (6, 3)      {$C_F: ABF$};
     \draw[->, thick] (1) -- (2);
     \draw[->, thick] (2) -- (3);
     \draw[->, thick] (2) -- (6);
     \draw[->, thick] (3) -- (5);
     \draw[->, thick] (6) -- (7);

\end{tikzpicture}
\caption{An RJT after implementing Line 10 for $n = A$.} \label{fig:rjt_example2}
\end{figure}

\begin{figure}[ht]
\centering 
\begin{tikzpicture}
    [cluster/.style={fill=white!80, draw, minimum size=2.5em, inner sep=2pt, rounded corners}]
     \node[cluster] (1) at (0, 2)      {$C_A: A$};
     \node[cluster] (2) at (2, 2)      {$C_B: AB$};
     \node[cluster] (3) at (4, 1)      {$C_C: ABC$};
     \node[cluster] (5) at (6, 1)      {$C_E: ABCE$};
     \node[cluster] (6) at (4, 3)      {$C_D: ABD$};
     \node[cluster] (7) at (6, 3)      {$C_F: ABF$};
     \draw[->, thick] (1) -- (2);
     \draw[->, thick] (2) -- (3);
     \draw[->, thick] (2) -- (6);
     \draw[->, thick] (3) -- (5);
     \draw[->, thick] (6) -- (7);

\end{tikzpicture}
\caption{An RJT after implementing Line 8 for $n = E$.} \label{fig:rjt_example3}
\end{figure}

\begin{figure}[ht]
\centering 
\begin{tikzpicture}
    [cluster/.style={fill=white!80, draw, minimum size=2.5em, inner sep=2pt, rounded corners}]
     \node[cluster] (1) at (0, 2)      {$C_A: A$};
     \node[cluster] (2) at (2, 2)      {$C_B: AB$};
     \node[cluster] (3) at (4, 2)      {$C_C: ABC$};
     \node[cluster] (5) at (6, 2)      {$C_E: ABCE$};
     \node[cluster] (6) at (2, 0)      {$C_D: ABD$};
     \node[cluster] (7) at (4, 0)      {$C_F: ABF$};
     \draw[->, thick] (1) -- (2);
     \draw[->, thick] (2) -- (3);
     \draw[->, thick] (5) -- (6);
     \draw[->, thick] (3) -- (5);
     \draw[->, thick] (6) -- (7);

\end{tikzpicture}
\caption{An RJT after implementing Line 9 for $n = E$.} \label{fig:rjt_example4}
\end{figure}

\begin{figure}[ht]
\centering 
\begin{tikzpicture}
    [cluster/.style={fill=white!80, draw, minimum size=2.5em, inner sep=2pt, rounded corners}]
     \node[cluster] (1) at (0, 2)      {$C_A: A$};
     \node[cluster] (2) at (2, 2)      {$C_B: AB$};
     \node[cluster] (3) at (4, 2)      {$C_C: ABC$};
     \node[cluster] (5) at (6, 2)      {$C_E: ABCE$};
     \node[cluster] (6) at (2, 0)      {$C_D: ABDE$};
     \node[cluster] (7) at (4, 0)      {$C_F: ABEF$};
     \draw[->, thick] (1) -- (2);
     \draw[->, thick] (2) -- (3);
     \draw[->, thick] (5) -- (6);
     \draw[->, thick] (3) -- (5);
     \draw[->, thick] (6) -- (7);

\end{tikzpicture}
\caption{An RJT after implementing Line 10 for $n = E$.} \label{fig:rjt_example5}
\end{figure}

\section{Problem size}
\label{subsec:problemsize}

From Definition 2.1, we can derive a relationship between the width of the tree and the size of the corresponding model. By definition, a tree with a width $k$ has a maximum cluster $C_n$ containing $k+1$ nodes. In a gradual RJT, the cluster $C_n$ includes exactly one node $n \in N$ not contained in its parent cluster $C_i$. Using the running intersection property, the other $k$ nodes in $C_n$ must also be in $C_i$. Assuming that all nodes $l \in C_i \cup C_n$ have at least two states $s_l$, this implies that there are $\prod_{l \in C_i \cup C_n} | S_n | \ge 2^k$ local consistency constraints (4) for the pair $(C_i, C_n) \in \mathscr{A}$ and the number of constraints in the model (2)-(8) is thus at least $\mathcal{O}(c^k)$, where $k$ is the width of the gradual RJT. This is in line with \citep{parmentier2020integer} pointing out that the RJT-based approach is only suited for problems with moderate rooted treewidth.

The width of the tree in Figure 4 (in the manuscript) is $N+1$, where $N$ is the number of treatment periods in the pig farm problem ($N = 3$ in the example), while the width of the original pig farm RJT in Figure 2 (manuscript) is only 2. Furthermore, we note that the rooted treewidth of a problem is defined as the size of the largest cluster minus one. In an RJT, we have $(I(n) \cup n) \subseteq C_n$ for all $n \in N$ and the treewidth is thus at least $\max_{n \in N}|I(n)|$. For the single-value-node pig farm problem, $|I(\bar{V})|=N+1$ and for the original pig farm problem, $|I(H_2)|=2$. Therefore, we conclude that there are no RJT representations for these problems with a smaller width than the ones presented in Figures 2 and 4 (manuscript).

These results imply that the optimization model for the original pig farm problem grows linearly with the number of stages, but the single-value-node pig farm grows exponentially with the number of decisions, suggesting possible computational challenges in larger problems.

\section{Proofs}
\label{sec:appendix}

\noindent\textbf{Proof of Proposition 3.2.} Consider a state combination $s_N \in S_N$. Let $g(s_N) = s_N$, where $g(s_N)_{\bar v} = (s_v)_{v \in N^{V}}$ and denote $\bar{s}_N = g(s_N)$ for compactness. Notice that the mapping $g$ is a bijection. The utility in graph $G$ is calculated as $\sum_{v \in N^V}u(s_v)$. In graph $\overline{G}$, the utility is $u(s_{\overline{v}}) = \sum_{v \in N^V}u(s_{v})$, which is equivalent to the utility evaluated from graph $G$. The probabilities are calculated as follows:

\begin{align*}
&\mathbb{P}_G(s_{N}) = \prod_{c \in N^{c} \cup N^{v}}\mathbb{P}_G(s_{c} \mid  s_{I(c)} )\prod_{d \in N^{D}}\delta(s_d \mid s_{I(d)}) &\\ 
&=\prod_{v \in N^{v}}\mathbb{P}_G(s_{v} \mid  s_{I(v)} ) \prod_{c \in N^{c}}\mathbb{P}_G(s_{c} \mid  s_{I(c)} )\prod_{d \in N^{D}}\delta(s_d \mid s_{I(d)}) & \\
&=\prod_{v \in N^{v}}\mathbb{P}_G(s_{v} \mid  s_{I(v)} )\prod_{c \in \overline{N}^{c}}\mathbb{P}_{\overline{G}}(\overline{s}_c \mid  \overline{s}_{I(c)} )\prod_{d \in \overline{N}^{D}}\delta(\overline{s}_d \mid \overline{s}_{I(d)}) &\\
&=\mathbb{P}_{\overline{G}}( \overline{s}_{\overline{v}} \mid  \overline{s}_{I(\overline{v})} )\prod_{c \in \overline{N}^{c}}\mathbb{P}_{\overline{G}}(\bar{s}_{c} \mid \bar{s}_{I(c)} )\prod_{d \in \overline{N}^{D}}\delta(\bar{s}_d \mid \bar{s}_{I(d)}) & \\
&=\mathbb{P}_{\overline{G}}(\bar{s}_{N}).&
\end{align*}
The second line of the proof follows directly from the first line by separating the chance nodes and value nodes in the first product. The third line of the proof follows directly from Lines 3-5 of Algorithm 1 and the definition of $g(s_N)$. Essentially, chance nodes and decision nodes (as well as their parents) remain the same in both diagrams and therefore, evaluating the product of decision strategy variables and conditional probabilities in chance nodes using the same states must be equivalent. To understand the fourth line of the proof, notice that $\overline{s}_{\overline{v}} = (s_{v_1})_{v \in N^V}$, which follows from the definition of $g(s_N)$. In essence, for each possible state combination of value nodes in $G$, there is a corresponding state for $\overline{v}$ in diagram $\overline{G}$. The equivalence $\mathbb{P}_{\overline{G}}( \overline{s}_{\overline{v}} \mid  \overline{s}_{I(\overline{v})} ) = \prod_{v \in N^{v}}\mathbb{P}_G(s_{v} \mid  s_{I(v)} )$ follows directly from Line 9 of Algorithm 1. The fifth line is trivial. Hence, the diagrams are equivalent, cf. Definition 3.1.
\\

\noindent\textbf{Proof of Proposition 3.3.} Assume that $M \subseteq N$. Algorithm 2 does not remove any clusters nor remove nodes from any cluster. In addition, if for any $a \in N$, $C_a \subsetneq C'_a$, then it must be that $n \preccurlyeq a, \forall n \in C'_a \setminus C_a$. Hence, Definition 2.1 $(b)$ and $(c)$ follow directly from the assumption. 

Next, property $(a)$ in Definition 2.1 (the running intersection property) is shown to hold for each line of the algorithm separately. Assume that for each line of the algorithm, an RJT that fulfils the running intersection property is given. In what follows, we will denote $C'_c$ as the root cluster of $c$ before implementing the line defined in Algorithm 2, and $C^{''}_c$ the cluster after implementing the line defined in Algorithm 2.

Lines 2, 3, 6 and 7 do not modify the RJT at all. For Line 1, Definition 2.1 $(a)$ holds trivially. Also for Line 10, the running intersection property follows directly from the assumption that the junction tree before Line 10 fulfils the running intersection property, since $C'_a \cap C'_b \subseteq C^{''}_a \cap C^{''}_b \subseteq   C^{''}_c = C'_c \cup \{n\}, \forall C'_c \in P(C'_a,C'_b), C'_a,C'_b \in P(C'_n,C'_m)$.

Line 8 adds $C'_e \cap C'_g$ to every cluster on path $P(C'_e,C'_n)$. Naturally, then for any two clusters of that path $C'_a,C'_b \in P(C'_e,C'_n)$, $C^{''}_a \cap C^{''}_b \subseteq \{C'_a \cap C'_b\} \cup \{C'_e \cap C'_g\} \subseteq C^{''}_c, \forall C^{''}_c \in P(C^{''}_a,C^{''}_b)$. The first subset sign follows from the assumption and definition of Line 8. The second subset sign follows from the fact, that $C'_a \cap C'_b \subseteq C^{''}_c$ by assumption and $C'_e \cap C'_g$ is added to $C^{''}_c$ in Line 8. Therefore, running intersection property holds after implementing Line 8.

For Line 9, consider $C'_a, C'_b \in P(C'_n,C'_m)$ such that $C'_e, C'_g \in P(C'_a,C'_b)$ for $e$ and $g$ as defined in Lines 6 and 7. By assumption, $C'_a \cap C'_b \in C'_c, \forall C'_c \in P(C'_a, C'_b)$, which includes clusters $C'_e$ and $C'_g$. Therefore, $C'_a \cap C'_b \in C'_e$ and $C'_a \cap C'_b \in C'_g \Rightarrow C'_a \cap C'_b \in C^{''}_c, \forall C^{''}_c \in P(C^{''}_e,C^{''}_g)$.

\section{$N$-monitoring example}

The $N$-monitoring problem from \citep{salo2022decision} represents a problem of distributed decision-making where $N$ decisions are made in parallel with no communication between the decision-makers. The node $L$ in Figure \ref{fig:n_monitoring_ID} represents a load on a structure, and nodes $R_i$ are reports of the load, based on which the corresponding fortification decisions $A_i$ are made. The probability of failure in node $F$ depends on the load and the fortification decisions, and the utility in $T$ comprises fortification costs and a reward if the structure does not fail.

\begin{figure}[ht]
    \centering
    \begin{tikzpicture}
    [decision/.style={draw, minimum size=2em, inner sep=2pt}, 
    chance/.style={circle, draw, minimum size=2em, inner sep=2pt},
    value/.style={diamond, draw, minimum size=2em, inner sep=2pt},
    scale=1.2]
    \node[chance]   (L) at (-0.25, 1.5)  {$L$};
    \node[chance]   (L1) at (0.75, 3)  {$R_1$};
    \node[chance]   (L2) at (0.75, 2.25)  {$R_2$};
    \node   (LN-1) at (0.75, 1)  {$\vdots$};
    \node[chance]   (LN) at (0.75, 0)  {$R_N$};
    \node[decision] (1) at (1.75, 3)  {$A_1$};
    \node[decision] (2) at (1.75, 2.25)  {$A_2$};
    \node (N-1) at (1.75, 1)  {$\vdots$};
    \node[decision] (N) at (1.75, 0)  {$A_N$};
    \node[chance]   (F) at (2.75, 1.5)  {$F$};
    \node[value]    (T) at (3.75, 1.5)  {$T$};     
    \draw[->, thick] (L) -- (L1);
    \draw[->, thick] (L) -- (L2);
    \draw[->, thick] (L) -- (LN);
    \draw[->, thick] (L) -- (F);
    \draw[->, thick] (L1) -- (1);
    \draw[->, thick] (L2) -- (2);
    \draw[->, thick] (LN) -- (N);
    \draw[->, thick] (1) -- (F);
    \draw[->, thick] (2) -- (F);
    \draw[->, thick] (N) -- (F);
    \draw[->, thick] (F) -- (T);
    \draw[->, thick] (1) -- (T);
    \draw[->, thick] (2) -- (T);
    \draw[->, thick] (N) -- (T);
\end{tikzpicture}
    \caption{An influence diagram representing the $N$-monitoring problem.}
    \label{fig:n_monitoring_ID}
\end{figure}

With topological order $L, R_1, A_1, ..., R_N, A_N, F, T$, the rooted junction tree corresponding to the diagram in Figure \ref{fig:n_monitoring_ID} is presented in Figure \ref{fig:n_monitoring_RJT}. The structure of parallel observations and decisions in the $N$-monitoring problem is very different compared to the partially observed Markov decision process (POMDP) structure of the pig farm problem.

\begin{figure}[ht]
    \centering
    \begin{tikzpicture}
    [cluster/.style={fill=white!80, draw, minimum size=2.5em, inner sep=2pt, rounded corners}]
     \node[cluster] (1) at (0, 0)           {$L$};
     \node[cluster] (2) at (2, 0)           {$L R_1$};
     \node[cluster] (3) at (4, 0)           {$L R_1 A_1$};
     \node[cluster] (4) at (0, -2)          {$L A_1 ... A_{k-1} R_k$};
     \node[cluster] (5) at (4, -2)          {$L A_1 ... A_{k-1} R_k A_k$};
     \node[cluster] (6) at (0, -4)        {$L A_1 ... A_{N-1} R_N A_N$};
     \node[cluster] (8) at (2.5, -4)          {$L A_1 ... A_N F$};
     \node[cluster] (9) at (4.5, -4)        {$A_1 ... A_N F T$};
     \draw[->, thick] (1) -- (2);
     \draw[->, thick] (2) -- (3);
     \draw[->, thick, dashed] (3) -- (4);
     \draw[->, thick] (4) -- (5);
     \draw[->, thick, dashed] (5) -- (6);
     \draw[->, thick] (6) -- (8);
     \draw[->, thick] (8) -- (9);
\end{tikzpicture}
    \caption{A rooted junction tree representing the $N$-monitoring problem.}
    \label{fig:n_monitoring_RJT}
\end{figure}

Any RJT model based on the tree in Figure \ref{fig:n_monitoring_RJT} has a treewidth of $N+2$. Decision programming formulations on the other hand scale much better with respect to $N$ compared to the pig farm problem, since only 2 nodes are added to the problem when $N$ is increased by 1, whereas in the pig farm problem, 3 nodes are added.

\begin{figure}[ht]
    \centering
    \includegraphics[width=0.45\textwidth]{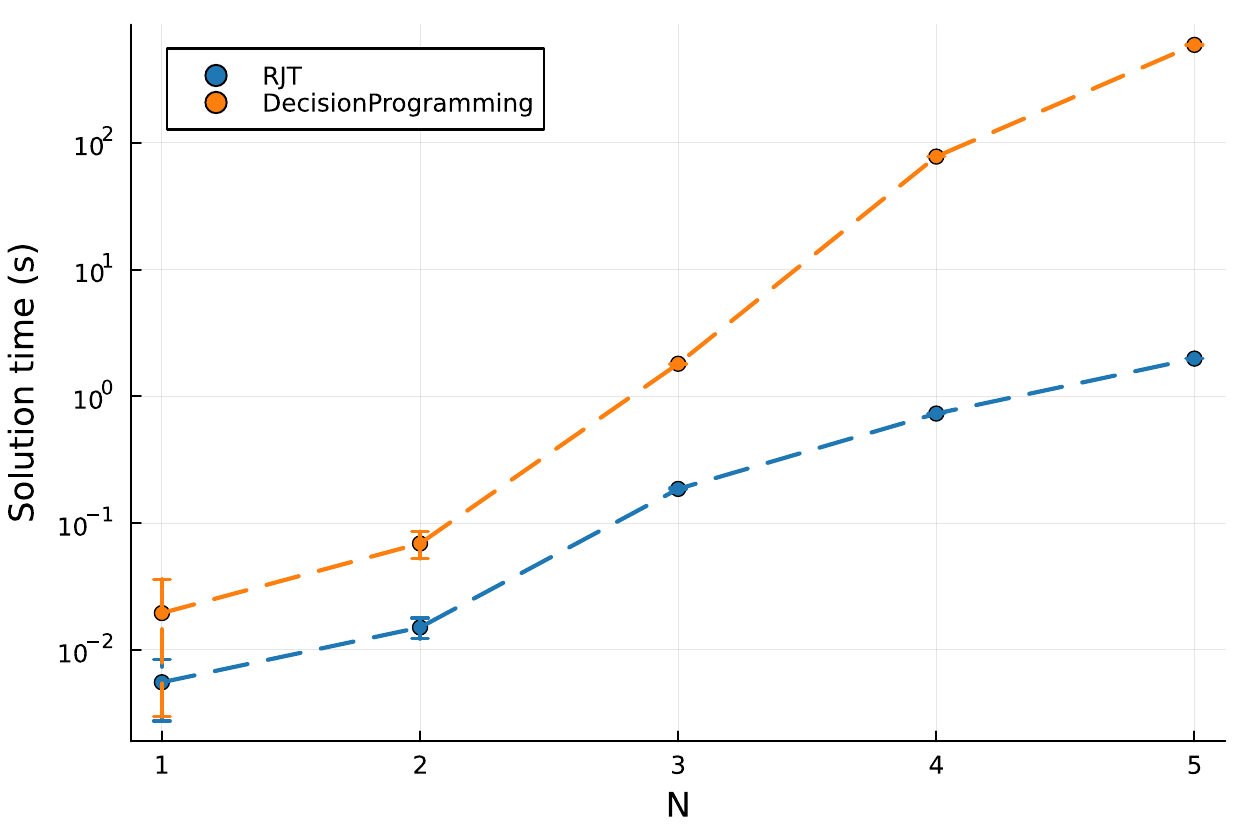}
    \caption{Mean solution times and standard deviations for 50 random instances of $N$-monitoring problem.}
    \label{fig:sol_times_DP_vsRJT_nmonitoring}
\end{figure}

Fifty randomly generated $N$-monitoring problem instances maximizing CVaR with probability level $1-\alpha = 0.85$ are solved with decision programming and the RJT model. This gives insights into the performance of the RJT model compared to decision programming for a problem that does not scale as well for the RJT model. The results are presented in Figure \ref{fig:sol_times_DP_vsRJT_nmonitoring}. The results show that the RJT formulation is more efficient in solving the $N$-monitoring problem. In larger instances, the difference in solution times is approximately 2 orders of magnitude. These results show that even though the structure of the problem would not be optimal for the RJT model, it still outperforms the corresponding decision programming model by a wide margin.

\bibliographystyle{elsarticle-harv} 
\bibliography{RJT_bibliography}